\documentclass[12pt]{amsart}
\usepackage{amsfonts}

\theoremstyle{definition}

\theoremstyle{remark}

\newcommand{\const}{\mathop{\rm const}\limits}

\let \vs = \vspace

\begin{document}

\begin{center}

{\bf Sharp Weight Hardy-Sobolev inequalities \\

\vspace{3mm}

for Grand Lebesgue Spaces.} \par

\vspace{4mm}

{\bf L. Sirota}\\

e-mail: sirota3@bezeqint.net \\

\vspace{4mm}

 \ Department of Mathematics and Statistics, Bar-Ilan University,\\
59200, Ramat Gan, Israel.\\

\vspace{4mm}

{\bf Abstract.} \\

\end{center}

{\it We give in this short paper the exact value for norms of two operators of Hardy-Sobolev type acting between two
weight Grand Lebesgue Space (GLS) based on the whole multidimensional Euclidean space.} \par

\vspace{3mm}


\vspace{3mm}


\vspace{3mm}

\section{Introduction.  Statement of problem.}

\vs{3mm}

 \  Let $  (X, |x|) $ be ordinary Euclidean space $  X = R^n $ with norm $ |x| = \sqrt{(x,x)}.  $ The standard Lebesgue-Riesz
$ L(p), \ p \ge 1 $ norm of a measurable function $ f: X \to R $ will be denoted

$$
| \ f \ |_p \stackrel{def}{=} \left[ \int_X | \ f(x) \ |^p \ dx  \right]^{1/p}.
$$

 \hspace{3mm} We recall for beginning  some famous and interest inequalities.  \par

\vspace{4mm}

{\bf A.  Hardy-Rellich inequality. } \\

\vspace{2mm}

 \ Introduce an (linear) operator (Hardy-Rellich operator)

$$
V[f](x) := \frac{f(x)}{|x|^2},
$$
then the Hardy-Rellich (-Sobolev) inequality may be written as follows, see \cite{Beckner1}

$$
| \ V[f] \ |_p \le K_{HR}(n,p) \ | \ \Delta f \ |_p, \eqno(1.1)
$$
 where $ \Delta $ is the Laplacian and

$$
K_{HR}(n,p) := \frac{p p'}{n(n - 2p)}, \hspace{4mm} p' := p/(p-1), \eqno(1.1a)
$$
and $ n \ge 3, \ 1 < p < n/2;   $ wherein this constant $ K_{HR}(n,p) $ is the best possible:

$$
\sup_{f: \ \Delta f  \ne 0} \frac{| \ V[f] \ |_p}{| \ \Delta f \ |_p} = K_{HR}(n,p). \eqno(1.2)
$$

 \ Obviously,  when $  p \to 1 + 0  $

$$
K_{HR}(n,p) \sim \frac{(p-1)^{-1}}{n(n-2)},
$$
and if $  p \to n/2 -0 $

$$
K_{HR}(n,p) \sim \frac{n}{4(n-2)} \cdot \frac{1}{n/2 -p}.
$$

 \ Recall that here  $ n \ge 3.$ \par

\vspace{4mm}

{\bf B.  Weight Sobolev's inequality. } \\

\vspace{3mm}

 \ Introduce an (linear) operator

 $$
 W[f](x) = |x|^{-\beta} \ f(x), \ x \in R^n, \ \beta = \const \in (0,n), \ n \ge 2, \ p \in (1, n/\beta);
 $$
then  \cite{Beckner1}

$$
| \ W[f](x) \ |_p  \le K_S(n, \beta,p) \cdot  | \ (- \Delta)^{\beta/2} \ f  \ |_p, \eqno(1.3)
$$
where the fractional degree of the anti-Laplace operator $ (- \Delta)^{\beta/2} $ is defined through the Fourier
transform and the best possible "constant"  $ K_S(n, \beta,p) $ has a form

$$
 K_S(n, \beta,p) = 2^{-\beta} \
 \frac{\Gamma\left[(\beta/(2p)) \cdot (n/\beta - p)\right] \cdot \Gamma \left[n/2 - n/(2p) \right]}
 {\Gamma \left[(n+\beta)/2 - n/(2p)  \right]  \Gamma \left[n/(2p) \right]}. \eqno(1.3a)
$$

\vspace{3mm}

 \ Note that as $ p \to 1 + 0  $

$$
K_S(n, \beta,p) \sim 2^{ 1 - \beta} \ \frac{\Gamma [ (n - \beta)/2 ] \ n^{-1} }{ \Gamma(\beta/2) \ \Gamma(n/2)} \cdot
(p-1)^{-1},
$$
and as $ p \to n/\beta - 0 $

$$
K_S(n, \beta,p)  \sim 2^{-\beta} \cdot \frac{\beta^2}{2n(n/\beta - p)} \cdot
\frac{\Gamma [ (n - \beta)/2  ]}{\Gamma(\beta/2) \ \Gamma(n/2)}.
$$

\vs{4mm}

 \ {\bf  Our aim is a generalization of the estimation (1.1) and (1.3) on the so - called
Grand Lebesgue Spaces $ GLS = GLS(\psi) = G(\psi), $ i.e. when } $ f(\cdot)
 \in G(\psi) \ $  {\bf and to show the accuracy of obtained estimations. } \par

\hfill $\Box$ \\
\bigskip

\section{Briefly about Grand Lebesgue Spaces.}

\vs{3mm}

 \  We recall briefly the definition and needed properties of these spaces.
  More details see in the works \cite{Fiorenza1}, \cite{Fiorenza2}, \cite{Ivaniec1},
   \cite{Ivaniec2}, \cite{Ostrovsky1}, \cite{Ostrovsky2}, \cite{Kozatchenko1},
  \cite{Jawerth1}, \cite{Karadzov1} etc. More about rearrangement invariant spaces
  see in the monographs \cite{Bennet1}, \cite{Krein1}. \par

\vspace{3mm}

 \ For $a$ and $b$ constants, $1 \le a < b \le \infty,$ let $\psi =
\psi(p),$ $p \in (a,b),$ be a continuous positive
function such that there exists a limits (finite or not)
$ \psi(a + 0)$ and $\psi(b-0),$  with conditions $ \inf_{p \in (a,b)} > 0 $ and
 $\min\{\psi(a+0), \psi(b-0)\}> 0.$  We will denote the set of all these functions
 as $ \Psi(a,b). $ \par

 \ The  Grand Lebesgue Space (in notation GLS) $  G(\psi; a,b) =
 G(\psi) $ is the space of all measurable
functions $ \ f: R^d \to R \ $ endowed with the norm

$$
||f||G(\psi) \stackrel{def}{=}\sup_{p \in (a,b)}
\left[ \frac{ |f|_p}{\psi(p)} \right], \eqno(2.1)
$$
if it is finite.\par
  \ In the article \cite{Ostrovsky2} there are many examples of these spaces.
 For instance, in the case when  $ 1 \le a < b < \infty, \beta, \gamma \ge 0 $ and

 $$
 \psi(p) = \psi(a,b; \beta, \gamma; p) = (p - a)^{-\beta} (b - p)^{-\gamma};
 $$
we will denote
the correspondent $ G(\psi) $ space by  $ G(a,b; \beta, \gamma);  $ it
is not trivial, non-reflexive, non-separable
etc.  In the case $ b = \infty $ we need to take $ \gamma < 0 $ and define

$$
\psi(p) = \psi(a,b; \beta, \gamma; p) = (p - a)^{-\beta}, p \in (a, h);
$$

$$
\psi(p) = \psi(a,b; \beta, \gamma; p) = p^{- \gamma} = p^{- |\gamma|}, \ p \ge h,
$$
where the value $ h $ is the unique  solution of a continuity equation

$$
(h - a)^{- \beta} = h^{ - \gamma }
$$
in the set  $ h \in (a, \infty). $ \par

 \ The  $ G(\psi) $ spaces over some measurable space $ (X, F, \mu) $
with condition $ \mu(X) = 1 $  (probabilistic case)
appeared in \cite{Kozatchenko1}.\par
 \ The GLS spaces are rearrangement invariant spaces and moreover interpolation spaces
between the spaces $ L_1(R^d) $ and $ L_{\infty}(R^d) $ under real interpolation
method \cite{Carro1}, \cite{Jawerth1}. \par
 \ It was proved also that in this case each $ G(\psi) $ space coincides
with the so-called {\it exponential Orlicz space,} up to norm equivalence. In others
quoted publications were investigated, for instance, their associate spaces, fundamental functions
$\phi(G(\psi; a,b);\delta),$ Fourier and {\it singular} integral operators acting in these spaces,
conditions for convergence and compactness, reflexivity and separability, martingales in these  spaces  etc.\par

\vs{3mm}

 \ {\bf Remark 2.1.} If we introduce the {\it discontinuous} function

$$
\psi_r(p) = 1, \ p = r; \hspace{3mm} \psi_r(p) = \infty, \ p \ne r, \ p,r \in (a,b)
$$
and define formally  $ C/\infty = 0, \ C = \const \in R^1, $ then  the norm
in the space $ G(\psi_r) $ coincides with the $ L_r $ norm:

$$
||f||G(\psi_r) = |f|_r.
$$

 \ Thus, the  Grand Lebesgue Spaces are the direct generalization of the
classical exponential Orlicz's spaces and Lebesgue spaces $ L_r. $ \par

\vs{3mm}

 \ The function $ \psi(\cdot) $ may be generated as follows. Let $ \xi = \xi(x)$
be some measurable function: $ \xi: X \to R $ such that $ \exists  (a,b):
1 \le a < b \le \infty, \ \forall p \in (a,b) \ |\xi|_p < \infty. $ Then we can
choose

$$
\psi(p) = \psi_{\xi}(p) = |\xi|_p. \eqno(2.2)
$$

 \ Analogously let $ \xi(t,\cdot) = \xi(t,x), t \in T, \ T $ is arbitrary set,
be some {\it family } $ F = \{ \xi(t, \cdot) \} $ of the measurable functions:
$ \forall t \in T  \ \xi(t,\cdot): X \to R $ such that
$$
 \exists  (a,b): 1 \le a < b \le \infty, \ \sup_{t \in T} \
|\xi(t, \cdot)|_p < \infty.
$$
Then we can choose

$$
\psi(p) = \psi_{F}(p) = \sup_{t \in T}|\xi(t,\cdot)|_p. \eqno(2.2a)
$$

 \ The function $ \psi_F(p) $ may be called as a {\it natural function} for the family $ F. $
This method was used in the probability theory, more exactly, in
the theory of random fields, see \cite{Ostrovsky1}. \par

 \ The GLS norm estimates, in particular, Orlicz norm estimates for
measurable functions, e.g., for random variables are used in the theory of
Partial Differential Equations \cite{Fiorenza1}, \cite{Ivaniec1}, theory of
probability in Banach spaces  \cite{Kozatchenko1},
\cite{Ostrovsky1}, in the modern non-parametrical statistics, for
example, in the so-called regression problem \cite{Ostrovsky1}.\par

\vspace{3mm}

 \ We use the symbols $C(X,Y),$ $C(p,q;\psi),$ etc., to denote positive
constants along with parameters they depend on, or at least
dependence on which is essential in our study. To distinguish
between two different constants depending on the same parameters
we will additionally enumerate them, like $C_1(X,Y)$ and
$C_2(X,Y).$ The relation $ g(\cdot) \asymp h(\cdot), \ p \in (A,B), $
where $ g = g(p), \ h = h(p), \ g,h: (A,B) \to R_+, $
denotes as usually

$$
0< \inf_{p\in (A,B)} h(p)/g(p) \le \sup_{p \in(A,B)}h(p)/g(p)<\infty.
$$

 \ The symbol $ \sim $ will denote usual equivalence in the limit
sense.\par

\hfill $\Box$ \\
\bigskip

\section{Main result: norm estimations for considered operators.}

\vspace{3mm}

 \ {\bf A.  Hardy-Rellich case. } \par

\vspace{3mm}

 \ Suppose  the (measurable) function $ \Delta f(\cdot): R^n \to R   $ belongs to some space
 $  G\psi(a,b),  $ i.e. $ || \Delta f||G\psi < \infty, $ where
  $ \psi(\cdot) \in \Psi(a,b), \ 1 \le a < b \le \infty.  $  For instance, the function $  \psi(\cdot)  $
 may be picked  as a natural function for $ \Delta f(\cdot), $ if there exists: $ \psi(p) = \psi_0(p),  $ where

$$
\psi_0(p) := | \ \Delta  f \ |_p, \hspace{5mm} p \in (a,b).
$$

 \ The set of all such a functions $ f:  \ f \in W^2 G\psi  = W^2 G\psi(R^n) $ equipped with the semi - norm

$$
|| \ f \ ||W^2 G\psi \stackrel{def}{=} || \ \Delta f \ ||G\psi
$$
 is named as usually Sobolev-Grand Lebesgue Space; it is a complete Banach space. \par

 \ Define the segment

$$
 I_{H R} = (p_0, p_1) \stackrel{def}{=} (1, n/2) \cap (a,b) \eqno(3.1)
$$
 and suppose  its non-triviality: $  I_{H R} \ne \emptyset  $ or equally  $  1 < p_0 < p_1 < n/2.  $ \par

 \ The most interesting and important case appears, by our opinion, when $ a = 1 $  and $  b = n/2; $ then

$$
  I_{H R} = (1, n/2).  \eqno(3.2)
$$

 \ Let us define a new $  \Psi \ - $ function

$$
\psi_V(p) \stackrel{def}{=} K_{HR}(p) \cdot \psi(p), \ \hspace{5mm} p \in  I_{H R}. \eqno(3.3)
$$

\vs{3mm}

 \ {\bf Theorem 3.A.}

\vs{3mm}

$$
|| \ V[f] \ ||G\psi_V \le 1 \cdot || \ f \ ||W^2 G\psi, \eqno(3.4)
$$
{\it where the constant "1" in (3.4) is the best possible in the case (3.2).} \par

\vs{3mm}

{\bf Proof.} Let $  f(\cdot) \in W^2 G\psi;  $  we can and will suppose without loss of generality

$$
|| \ f \ ||W^2 G\psi = 1.
$$

  \ It follows immediately from the direct of the Grand Lebesgue Spaces norm

$$
| \ \Delta f \ |_p \le \psi(p), \hspace{5mm} p \in (a,b).
$$

 \ We use the inequality (1.1):

$$
| \ V[f] \ |_p \le K_{HR}(n,p) \ | \ \Delta f \ |_p \le K_{HR}(n,p) \cdot \psi(p) =
$$

$$
\psi_V(p) =   \psi_V(p) \cdot  || \ f \ ||W^2 G\psi(p),
$$
or equally

$$
|| \ V[f] \ ||G\psi_V \le || \ f \ ||W^2 G\psi.
$$

 \ The sharpness  of the constant $  "1" $ follows immediately from  one of the results of  the preprint
\cite{Ostrovsky3}. \par

\vs{4mm}

 \ {\bf B.  Sobolev's case. } \par

\vs{4mm}

 \ Suppose now that  the (measurable) function $  f(\cdot): R^n \to R   $  is such that $  (-\Delta)^{\beta/2} f  $
 belongs to some space $  G\psi(a,b),  $ i.e. $ || (-\Delta)^{\beta/2} f||G\psi < \infty, $ where
  $ \psi(\cdot) \in \Psi(a,b), \ 1 \le a < b \le \infty.  $  For instance, the function $  \psi(\cdot)  $
 may be picked  as a natural function for $ (- \Delta)^{\beta/2} f(\cdot), $ if there exists: $ \psi(p) = \psi_1(p),  $ where

$$
\psi_1(p) := | \ (-\Delta)^{\beta/2}  f \ |_p, \hspace{5mm} p \in (a,b).
$$

 \ The set of all such a functions $ f: \ f \in W^{(\beta)} G\psi  = W^{(\beta)} G\psi(R^n) $ equipped with the semi - norm

$$
|| \ f \ ||  W^{(\beta)} G\psi \stackrel{def}{=} || \ (-\Delta f)^{\beta/2} \ ||G\psi
$$
 is named again as usually Sobolev-Grand Lebesgue Space; it is also the complete Banach space. \par

 \ Define the segment

$$
 J_{W} = (p_2, p_3) \stackrel{def}{=} (1, n/\beta) \cap (a,b) \eqno(3.5)
$$
 and suppose  its non-triviality: $  J_{W} \ne \emptyset  $ or equally  $  1 < p_2 < p_3 < n/\beta.  $ \par

 \ The most interesting and important case appears, by our opinion, when $ a = 1 $  and $  b = n/\beta; $ then

$$
  J_{H R} = (1, n/\beta).  \eqno(3.6)
$$

 \ Let us define a new $  \Psi \ - $ function

$$
\psi_W(p) \stackrel{def}{=} K_{W}(p) \cdot \psi(p), \ \hspace{5mm} p \in  J_{H R}. \eqno(3.7)
$$

\vs{3mm}

 \ Analogously to the theorem 3.A may be proved the following result. \par

\vspace{3mm}

 \ {\bf Theorem 3.B.}

\vs{3mm}

$$
|| \ W[f] \ ||G\psi_W \le 1 \cdot || \ f \ ||W^{(\beta)} G\psi, \eqno(3.8)
$$
{\it where the constant "1" in (3.8) is the best possible in the case (3.6).} \par

\hfill $\Box$ \\
\bigskip

 \vspace{3mm}

\section{Simplifications.}

\vs{3mm}

 \ Notice that we know

$$
K_{HR}(p) \asymp \frac{1}{(p-1)(n/2 - p)}, \hspace{4mm} p   \in (1, n/2).
$$
  \ More precisely,

$$
\frac{C_1(n)}{(p-1)(n/2 - p)} \le  K_{HR}(p) \le \frac{C_2(n)}{(p-1)(n/2 - p)}, \hspace{4mm} p   \in (1, n/2),
$$
 where as before $  n \ge 3 $  and

$$
0 < \inf_{n \ge 3}  C_1(n) < \sup_{n \ge 3} C_2(n) < \infty.
$$

 \ Analogously

$$
\frac{C_3(n,\beta)}{(p-1)(n/\beta - p)} \le  K_{S}(p) \le \frac{C_4(n,\beta)}{(p-1)(n/\beta - p)}, \hspace{4mm} p   \in (1, n/\beta),
$$
 where  $  n \ge 2 $  and

$$
0 < \inf_{n \ge 2}  C_3(n,\beta) < \sup_{n \ge 2} C_4(n,\beta) < \infty.
$$

 \ Therefore, one can in the relations  (3.3) and (3.7) the coefficients $ K_{HR}(p), \  K_{S}(p) $ replace to the more simple
expressions. The assertions of theorem (3.A) and (3.B) remain true up to finite multiplicative constants. \par

\hfill $\Box$ \\
\bigskip

 \section{Concluding remarks}

\vspace{3mm}

 \ Of course, described here method may be applied to the more wide  class of operators, non necessary to be linear,
differential or integral, for which are known $ L(p)-L(q) $  estimates.  Many examples of these operators one can found
in articles and books \cite{Beckner101},  \cite{Hardy1}, \cite{Karadzov1}, \cite{Liflyand1}, \cite{Maz'ya1}, \cite{Mitrinovich1},
\cite{Okikiolu1}, \cite{Ostrovsky3}-\cite{Ostrovsky4}, \cite{Perez1}, \cite{Stein1}, \cite{Talenti1} etc.

\end{document}